\documentclass{amsart}
\usepackage{amssymb,amsfonts,amsmath}
\usepackage[all]{xy}
\usepackage[shortlabels]{enumitem}
\usepackage{hyperref, color}

\usepackage[pagewise]{lineno}
		
\newtheorem{theorem}{Theorem}[section]

\newtheorem{corollary}{Corollary}[section]
\theoremstyle{definition}

\theoremstyle{lemma}
\newtheorem{lemma}{Lemma}[section]

\theoremstyle{remark}
\newtheorem{remark}{Remark}
\theoremstyle{claim}
\newtheorem{claim}{Claim}[section]

\theoremstyle{theoremA}

\theoremstyle{theoremB}
  
\newcommand{\real}{\mathbb{R}}

\newcommand{\hy}{\mathbb{H}}
\newcommand{\Sp}{\mathbb{S}}

\newcommand{\sech}{\mathrm{sech}}

\newcommand{\si}{\sigma}
\newcommand{\te}{\theta}

\newcommand{\al}{\alpha}

\newcommand{\ep}{\epsilon}
\newcommand{\ga}{\gamma}
\newcommand{\Ga}{\Gamma}
\newcommand{\be}{\beta}

\newcommand{\De}{\Delta}
\newcommand{\de}{\delta}

\newcommand{\ti}{\tilde}

\newcommand{\lan}{\left\langle}
\newcommand{\ran}{\right\rangle}
\newcommand{\p}{\partial}

\newcommand{\m}{\mathcal}

\newcommand{\vp}{\varphi}
\begin{document}


\title[Multiple warped product metrics in quadrics]{Global isometric embeddings of multiple warped product metrics into quadrics}
\thanks{The authors were partially supported by CNPq-Brazil}

\author{H. Mirandola}

\author{F. Vit\'orio}

\date{}

\maketitle
\thispagestyle{empty}
\begin{abstract} In this paper, we construct smooth isometric embeddings of multiple warped product manifolds  in quadrics of semi-Euclidean spaces. Our main theorem generalizes previous results as given by Blanu\v sa, Rozendorn, Henke and Azov. 
\end{abstract}
\section{Introduction}
By fundamental works of Nash \cite{na}, Gromov and Rokhlin \cite{gr} and G\"unther \cite{gu} we know that every $n$-dimensional smooth Riemannian manifold admits a smooth isometric embedding \footnote{We recall that an injective immersion is an embedding if it is a homeomorphism onto its image, by considering the image with the induced topology.} in an $N$-dimensional Euclidean space $\real^N$, for some $N\leq c(n)=\max\{n(n+5)/2, n(n+3)/2+5\}$. The estimate $c(n)$ was given by G\"unther \cite{gu}; Nash's and Gromov-Rokhlin's estimates are larger than this upper bound. Since then, the problem of finding the lowest possible codimension is one of the major open problems in the theory of  isometric immersion. For books and surveys about this subject see Gromov and Rokhlin \cite{gr}, Jacobowitz \cite{jac}, Poznyak and Sokolov \cite{poso}, Aminov \cite{am}, Dajczer \cite{da}, Borisenko \cite{bo} and Han and Hong \cite{hh}.  On the other hand, since the results of Nash, Gromov and Rokhlin and G\"unther follow as a consequence of existence theorems for certain PDE's, it is also an interesting problem to give the explicit construction of isometric immersions of a given Riemannian metric  $M^n$ in $\real^m$, mainly if the attained codimension is strictly less than $c(n)-n$. This is the point of view of the present paper.

Blanu\v sa \cite{bla, bla2} gave a method to construct injective smooth isometric immersions of the hyperbolic plane $\hy^2$ in $\real^6$ and in the standard spherical space $\Sp^8$. Poznyak \cite{po} wrote about Blanu\v sa's surface: {\it{``There is no doubt that this result is one of the most elegant in the theory of immersion of two-dimensional manifolds in Euclidean space"}}. Blanu\v sa also constructed injective
isometric immersions of $\hy^n$ in $\real^{6n-5}$ and of an infinite M\"obius band with hyperbolic metric in $\real^8$ and in $\Sp^{10}$. 
His method was used and modified 
in further works: (i) Rozendorn \cite{rozen} constructed non-injective smooth isometric immersions of the plane $\real^2$ with the warped product metric  of the form $d\si^2=dt^2+f(t)^2dx^2$ in $\real^5$ (this class of surfaces includes $\hy^2$ with the metric $dt^2+e^{2t}dx^2$). Note that a celebrated theorem of Hilbert \cite{hilb} states that $\hy^2$ cannot be isometrically immersed in $\real^3$. However, the existence of an isometric immersion of $\hy^2$ in $\real^4$ or even an injective isometric immersion of $\hy^2$ in $\real^5$ is still an open problem (a partial answer to the first problem was given by Sabitov \cite{sab}). (ii) Henke \cite{henke1, henke3} exhibited isometric immersions of $\hy^n$ in $\real^{4n-3}$ and in  $\Sp^{4n-3}$.  (iii) Henke and Nettekoven \cite{henke2} showed that $\hy^n$ can be isometrically embedded in $\real^{6n-6}$ whose image is the graph of a smooth map $g:\real^n\to \real^{5n-6}$. (iv) Azov \cite{azov} considered the space $\real^n=\real\times \real^{n-1}$ with one of the following metrics: $d\si^2 = dt^2+ f(t)^2\sum_{j=1}^{n-1}dx_j^2$ or $d\si^2 = g(x_1)^2\sum_{j=1}^n dx_j^2$ and constructed isometric immersions in $\real^{4n-3}$ and $\Sp^{4n-3}$. He also announced in \cite{azov2} the construction of isometric immersions of these classes of metrics in $\real^{4n-4}$ and $\Sp^{4n-4}$, if $n>2$.

In this paper we deal with product manifolds $M^n=I\times \real^{n-1}$, where $I$ is an open interval, endowed with a multiple warped product metric of the form:
\begin{eqnarray}\label{multiplewarped}
d\si^2 = \,\rho(t)^2\,dt^2 +\eta_1(t)^2 dx_1^2+\ldots +\eta_{n-1}(t)^2 dx_{n-1}^2
\end{eqnarray}
where $\rho(t),\eta_j(t)$, with $t\in I$ and  $j=1,\ldots,n-1$, are positive smooth functions and $dx_1,\ldots, dx_{n-1}$ are the canonical coframes of $\real^{n-1}$. This class of metrics includes both Azov's metrics. We will modify  Blanu\v sa's method to exhibit isometric immersions and, mainly, embeddings of this class of metrics in quadrics of semi-Euclidean spaces. It is worth to mention that, in general, the immersions obtained by Rozendorn,  Henke and Azov are not injective.  Based on this, we consider such embeddings the main contribution of the present work. There exists a wide literature about aspects of rigidity and nonimmersibility of these spaces (see for instance N\"olker \cite{no}, Chen \cite{ch}, Florit \cite{flo}, Dajczer and Tojeiro \cite{dato} and references therein).  

We recall that the semi-Euclidean space $\real^n_a$, with $a\in \{0,\ldots,n\}$, is simply the space $\real^n$  with the inner product of signature $(a,n-a)$ given by
\begin{equation}\label{minkowiski}
\lan\,,\ran =\,-\,dx_1^2-\ldots - dx_a^2 + dx_{a+1}^2 + \ldots + dx_n^2,
\end{equation}
where $dx_j$, with $j=1,\ldots, n$, denote the canonical coframes of $\real^n$. 
For a given $c>0$, let $S^n_a(c)$ and $H^n_a(-c)$ be the following quadratic hypersurfaces (or, simply, quadrics):
\begin{eqnarray*}
&&\Sp^n_a(c) =\left\{x\in \real^{n+1}_a \bigm| \lan x,x\ran={1}/{c}\right\};\\ 
&&\hy^n_a(-c) = \left\{x\in \real^{n+1}_{a+1} \bigm| \lan x,x\ran=-{1}/{c}\right\}.
\end{eqnarray*}
Both hypersurfaces are semi-Riemannian manifolds with signature $(a,n-a)$ and constant curvatures $c$ and $-c$, respectively. If $a=0$, then $\Sp^n_0(c)=\Sp^n(c)$ is the standard sphere and $\hy^n_0(-c)=\hy^n(-c)$ is the hyperbolic space. If $a=1$, the semi-Riemannian universal covering spaces of $\Sp^n_1(1)$ and $\hy^n_1(-1)$ are called de Sitter $dS^n$ and anti-de Sitter $AdS^n$ spaces, respectively.

Our main result is
\begin{theorem}\label{mth} Let $M^n=I\times \real^{n-1}$ be as given in (\ref{multiplewarped}). 
Then, for all $c>0$ and $a\in \{0,\ldots, n-1\}$, the manifold $M^n$ admits:
\begin{enumerate}[(i)]
\item\label{euc-4n-2l-3} isometric immersions in $\real^{4n-3-2a}_{a}$, $\Sp^{4n-3-2a}_a(c)$, and $\hy^{4n-3-a}_a(-c)$;
\item\label{imb-8n-6a-7} isometric embeddings in $\real^{8n-7-6a}_a$, $\Sp^{8n-5-6a}_a(c)$ and $\hy^{8n-7-5a}_a(-c)$.
\end{enumerate}
Moreover, all immersions and embeddings above are  smooth and given explicitly. 
\end{theorem}
In Remark \ref{non-inj} (see Section \ref{sec-proof-mth}), we observe that all immersions referred in Item {\it \ref{euc-4n-2l-3}} of Theorem \ref{mth} are not injective, provided that $a<n-1$. 

Based on Theorem \ref{mth}, it is natural to ask if every $n$-dimensional Riemannian manifold $M^n$ can be isometrically immersed in a semi-Euclidean $\real^N_a$ with $a>0$ and $N$ strictly less than the dimension $c(n)$ obtained by G\"unther \cite{gu}.

As an application of Theorem \ref{mth} we generalize  Rozendorn's surfaces \cite{rozen}. We have the following.
\begin{corollary} Let $M^2=I\times \real$ be a warped product surface as given in (\ref{multiplewarped}). Then, for all $c>0$, the surface $M^2$ admits:
\begin{enumerate} [(i)]
\item non-injective isometric immersions in $\real^5$, $\hy^5(-c)$, $\Sp^5(c)$ and $dS^3(c)$;
\item isometric embeddings in $\real^9$, $\hy^9(-c)$, $\Sp^{11}(c)$, $\real^3_1$, $AdS^4(-c)$ and $dS^5(c)$.
\end{enumerate}
Moreover, all immersions and embeddings above are smooth and given explicitly.
\end{corollary} 

The space $Sol_3$ is a simply connected homogeneous $3$-dimensional space whose isometry group has dimension $3$. It is one of the eight models of the Thurston geometry and it can be viewed as $\real^3$ with the metric $ds^2=dt^2 + e^{2t}dx^2+e^{-2t}dy^2.$ It follows directly from Theorem \ref{mth} the following

\begin{corollary}\label{two-dimensional} For all $c>0$, the space $Sol_3$ admits:
\begin{enumerate}[(i)]
\item non-injective isometric immersions in $\real^9$, $\hy^9(-c)$, $\Sp^9(c)$, $\real^7_1$, $dS^7(c)$, $AdS^8(-c)$;
\item isometric embeddings in $\real^{17}$, $\hy^{17}(-c)$, $\Sp^{19}(c)$, $\real^{11}_1$, $AdS^{12}(-c)$ and $dS^{13}(c)$.
\end{enumerate} 
Moreover, all immersions and embeddings above are smooth and given explicitly.
\end{corollary}

Let $f_l:M_l\to \real^{n_l}$, with $l=1,\ldots,k$, be smooth isometric immersions of the manifold $(M_l,g_l)$ in $\real^{n_l}$. Let $I$ be an open interval and $\rho(t),\eta_l(t)$, with $t\in I$ and $l=1,\ldots,k$,  positive smooth functions. It is simple to show that the product manifold $M=I\times M_1\times\ldots\times M_k$ with the warped product metric
\begin{equation}\label{warper metr riem}
g=\rho(t)^2 dt^2 + \eta_1(t)^2 g_1+\ldots+\eta_k(t)^2 g_k
\end{equation}
can be isometrically immersed in $I\times \real^{n_1}\times\ldots\times\real^{n_k}$ with the metric
\begin{equation*}
d\si^2 = \rho(t)^2 dt^2 + \eta_1(t)^2 \de_1+\ldots+\eta_k(t)^2\de_k,
\end{equation*}
where $\de_l$ denotes the Euclidean metric of $\real^{n_l}$. Thus it follows as a consequence of Theorem \ref{mth} the following result.

\begin{corollary} With the notations being as above, we consider $n=n_1+\ldots+n_k$ and $a\in \{0,\ldots,n\}$. For all $c>0$, the manifold $M$ admits:
\begin{enumerate}[(i)]
\item isometric immersions in $\real^{4n+1-2a}_{a}$, $\Sp^{4n+1-2a}_{a}(c)$ and $\hy^{4n+1-a}_{a}(-c)$;
\item isometric embeddings in $\real^{8n+1-6a}_{a}$, $\Sp^{8n+3-6a}_{a}(c)$ and $\hy^{8n+1-5a}_{a}(-c)$, provided that each $f_l$ is an embedding.
\end{enumerate}
\end{corollary}

We would like to thank the referee for carefully reading the first version of this manuscript, pointing out mistakes which helped us to improve the manuscript.

\section{Preliminaries.}\label{method}

We recall Blanu\v sa's functions $\hat\psi_1,\hat\psi_2:\real\to \real$  defined by
\begin{equation*}
\hat\psi_1(u)=\sqrt{\frac{1}{A}\int_0^{u+1}\xi(\tau)d\tau} \ \ \mbox{ and } \ \ \hat\psi_2(u)=\sqrt{\frac{1}{A}\int_0^{u}\xi(\tau)d\tau},
\end{equation*}
where $A=\int_0^{1}\xi(\tau)d\tau$ and $\xi(u)=\sin(\pi u)e^{\frac{-1}{(\sin(\pi u))^2}}$, if $u\in \real\setminus \Bbb Z$, and $\xi(u)=0$, if $u\in \Bbb Z$. 
Blanu\v sa proved in \cite{bla} that these functions are smooth, non-negative and satisfy: 
\begin{enumerate}[(a)]
\item\label{bla-it1} $\hat\psi_j$ is periodic with period 2, for all $j=1,2$; 
\item\label{bla-it2} $\hat\psi_1^2+\hat\psi_2^2=1$, everywhere; 
\item\label{bla-it3} all the derivatives $\hat\psi_1^{(k)}(2l+1)=\hat\psi_2^{(k)}(2l)=0$, for all $l\in \Bbb Z$.
\end{enumerate} 

The next two lemmas will be useful to prove Theorem \ref{mth}. They are simple consequences of Items \ref{bla-it1}, \ref{bla-it2} and \ref{bla-it3} above. To state them, let $I$ be an open interval and $\ga:I\to \real$ an increasing smooth diffeomorphism. Consider the sequence $\mathbf{t}_k=\ga^{-1}(k)$, with $k\in \Bbb Z$. The first lemma says the following.

\begin{lemma}\label{bla-fun-I} The functions $\psi_j=\hat\psi_j\circ \ga:I\to \real$ are smooth, non-negative and satisfy the following properties:
\begin{equation}\label{bla-fun}
\left\{
\begin{array}{l}
\psi_1(t)^2+\psi_2(t)^2=1, \mbox{ everywhere in } I;\\ 
\psi_j\big((\ga)^{-1}(u)\big)=\psi_j\big((\ga)^{-1}(u+2)\big), \mbox{ for all } u\in \real \mbox{ and } j=1,2;\\
\psi_1^{(k)}(\mathbf{t}_{2l+1})=\psi_2^{(k)}(\mathbf{t}_{2l})=0, \mbox{ for all }k\geq 0 \mbox{ and }  \mbox{ integers } l.
\end{array}\right.
\end{equation}
\end{lemma}

Let $S_1,S_2:I\to (0,\infty)$ be any positive step functions satisfying 
\begin{equation}\label{stp-fun}\left\{
\begin{array}{l}
S_1 \mbox{ is constant on each interval } [\mathbf{t}_{2l+1},\mathbf{t}_{2l+3});\\
S_2 \mbox{ is constant on each interval } [\mathbf{t}_{2l}, \mathbf{t}_{2l+2});
\end{array}\right.
\end{equation}
for each integer $l$.

The second lemma  follows easily from Lemma \ref{bla-fun-I}.

\begin{lemma}\label{lemmaderiv}
For any $\eta\in C^\infty(I)$,  the functions 
$\frac{\eta(t)\psi_j(t)}{S_j(t)}$ with $t\in I$ and $j=1,2,$
are smooth and their derivatives satisfy
\begin{equation*}
\frac{d^{\,k}}{dt^k}\left(\frac{\eta(t)\psi_j(t)}{S_j(t)}\right)=\frac{\frac{d^{\,k}}{dt^k}\Big(\eta(t)\psi_j(t)\Big)}{S_j(t)},
\end{equation*}
for all integers $k\geq 0$.
\end{lemma}

\section{Proof of Theorem \ref{mth}}\label{sec-proof-mth}
First we consider the map $\eta(t)=(\eta_1(t),\ldots,\eta_{n-1}(t))$, with $t\in I$, where each function $\eta_j$ is being as in (\ref{multiplewarped}). 
Consider the map $h:\real\to \real^2_1$  given by $h(u)=(\cosh(u),\sinh(u))$. Consider also the map $\vp=(\vp_1,\vp_2):I\times \real \to \real^4$  where each map $\vp_j:I\times \real\to \real^2$, with $j=1,2$, is given by 
\begin{equation}\label{map-euc}
\vp_j(t,u)=\frac{\psi_j(t)}{S_j(t)}\Big(\cos(S_j(t)u),\sin(S_j(t)u)\Big).
\end{equation}
The map $\vp$ is introduced in \cite{henke2} for the case that $I=\real$ and $\ga$ is the identity function. 
By using Lemma \ref{lemmaderiv}, we obtain  
\begin{eqnarray}\label{prelim}
&&\frac{\p(\eta_k(t) h(u))}{\p t} = \eta_k'(t)\big(\cosh(u),\sinh(u)\big);\\
&&\frac{\p(\eta_k(t) h(u))}{\p u} = \eta_k(t)\big(\sinh(u),\cosh(u)\big); \nonumber\\
&&\frac{\p (\eta_{k}(t)\vp_j(t,u))}{\p t} = \frac{(\eta_{k}(t)\psi_j(t))'}{S_j(t)}\Big(\cos(S_j(t)u),\sin(S_j(t)u)\Big); \nonumber\\
&&\frac{\p (\eta_{k}(t)\vp_j(t,u))}{\p u} = \eta_{k}(t)\psi_j(t)\Big(-\sin(S_j(t)u)\,,\,\cos(S_j(t)u)\Big);\nonumber
\end{eqnarray}
for all $j=1,2$, $t\in I$, $u\in \real$ and $k=1,\ldots,n-1$.

Now, set $a\in \{0,\ldots, n-1\}$ and let $b=n-1-a$. First we consider $b>0$. We will see that the case $b=0$ is easier. 
We write the semi-Euclidean space $\real^{4n-4-2a}_a=\real^{2a+4b}_{a}$ isometrically as the following form
\begin{equation*}
\real^{2a+4b}_{a}=(\real^{2}_1)^a\times \real^{4b} = \underbrace{\real^2_1 \times \ldots \times \real^2_1}_{\textrm{$a$  times }} \, \times \, \real^{4b}.
\end{equation*}
We denote by $x=(x_1,\ldots,x_{a+b})$ the coordinates of $\real^{n-1}=\real^{a+b}$. Let $P_1:\real^{a+b}\to \real^a$ and $P_2:\real^{a+b}\to \real^b$ be the standard orthogonal projections 
\begin{eqnarray*}
&&\ti x=P_1(x_1,\ldots,x_{a+b}) = (x_1,\ldots,x_a)\\ &&\bar x= P_2(x_1,\ldots,x_{a+b})=(x_{a+1},\ldots,x_{a+b}).
\end{eqnarray*}

Consider the maps 
\begin{equation}\label{ti-bar}
\ti\eta(t)=P_1(\eta(t)) \,\mbox{ and }\,\bar \eta(t)=P_2(\eta(t)),
\end{equation}
with $t\in I$ and let $\ti\eta\star h:I\times \real^{n-1}\to (\real^2_1)^a$ and $\bar\eta\star\vp:I\times \real^{n-1}\to \real^{4b}$ be the maps given by
\begin{eqnarray}
&&(\ti\eta\star h)(t,x_1,\ldots,x_{a+b}) = \big(\eta_1(t)h(x_1),\ldots, \eta_a(t) h(x_a)\big)\in (\real^2_1)^a;\label{starfunction}\\
&&(\bar\eta\star\vp)(t,x_{1},\ldots,x_{a+b}) = \big(\eta_{a+1}(t)\vp(t,x_{a+1}), \ldots, \eta_{a+b}(t)\vp(t,x_{a+b})\big)\in\real^{4b}\nonumber.
\end{eqnarray}
Since $h(u)\in \real^2_1$ and $\vp(t,u)\in \real^4$, by using (\ref{prelim}), the pull-back symmetric tensors by $\ti\eta\star h$ and $\bar\eta\star\vp$ must satisfy
\begin{eqnarray}\label{etah} 
(\ti\eta\star h)^*(\lan\,,\ran)&=& -|\ti\eta'(t)|^2\,dt^2 + \eta_1(t)^2 dx_1^2+\ldots+ \eta_a(t)^2 dx_a^2;\\
(\bar\eta\star\vp)^*(\lan\,,\ran) &=& 
\ep(t)^2\,dt^2 + \eta_{a+1}(t)^2 dx_{a+1}^2 + \ldots + \eta_{a+b}(t)^2 dx_{a+b}^2, \nonumber
\end{eqnarray}
where $|\ti\eta'(t)|^2=\eta_1'(t)^2+\ldots+\eta_a'(t)^2$ and $\ep:I\to [0,\infty)$ is  given by
\begin{equation}\label{ep_r}
\ep(t)^2= 
\sum_{r=a+1}^{a+b}\left[\frac{\left(\left(\eta_{r}(t)\psi_1(t)\right)'\right)^2}{S_1(t)^2} + \frac{\left(\left(\eta_{r}(t)\psi_2(t)\right)'\right)^2}{S_2(t)^2}\right].
\end{equation}

For the step functions $S_1,S_2:I\to (0,\infty)$ defined as in (\ref{stp-fun}), we can choose the steps $S_1|_{[\mathbf{t}_{2l+1}, \mathbf{t}_{2l+3})}$ and $S_2|_{[\mathbf{t}_{2l},\mathbf{t}_{2l+2})}$, with integer $l$, sufficiently large so that, for all $r=a+1,\ldots,a+b$, it holds 
\begin{eqnarray}\label{cond-ep}
\big(\left(\eta_{r}(t)\psi_j(t)\right)'\big)^2 <\frac{1}{4b}S_j(t)^2 \,\rho(t)^2,
\end{eqnarray}
for all $t\in I$ and $j=1,2$.
We obtain $\rho(t)^2-\ep(t)^2\geq \rho(t)^2-2\ep(t)^2> 0,$ for all $t\in I$. 

Let $f:I\times \real^{n-1}\to \real^{4n-3-2a}_a=\real\times(\real^2_1)^a\times \real^{4b}$ be the map  
\begin{equation*}
f(t,x)=
\left(\int_{t_0}^t \sqrt{\rho(\tau)^2+|\ti\eta'(\tau)|^2-\ep(\tau)^2}d\tau\,,\, \ti\eta\star h\,(t,x)\,,\,\bar\eta\star \vp\,(t,x)\right).
\end{equation*}
If $b=0$, we define $f(t,x)$ simply by omitting $\ep(t)$ and $\bar\eta\star \vp(t,x)$ in the expression of $f(t,x)$ above.
By using (\ref{etah}), we obtain  
\begin{eqnarray*}
f^*(\lan\,,\ran) &=& (\rho(t)^2+|\ti\eta'(t)|^2-\ep(t)^2) dt^2 + (\ti\eta\star h)^*(\lan\,,\ran) + (\bar\eta\star \vp)^*(\lan\,,\ran) \\ &=&\rho(t)^2 dt^2 + \eta_1(t)^2dx_1^2+\ldots + \eta_{a+b}(t)^2 dx_{a+b}^2\\
&=& d\si^2.
\end{eqnarray*}
This implies that $f:M^n\to \real^{4n-3-2a}_a$ is a smooth isometric immersion. 

We fix $c>0$. First we assume $b>0$. We choose the step functions $S_1,S_2$ sufficiently large so that (\ref{cond-ep}) holds. Let $\al:I\to [0,\infty)$ be the function given by
\begin{equation}\label{def-a}
\al(t) = \sum_{r=a+1}^{a+b}\left(\frac{\eta_r(t)^2\psi_1(t)^2}{S_1(t)^2}+\frac{\eta_r(t)^2\psi_2(t)^2}{S_2(t)^2}\right).
\end{equation}
Note that $\al(t)=\lan \bar\eta\star\vp\,(t,x), \bar\eta\star\vp\, (t,x)\ran$. 

Let $f_h:I\times \real^{n-1}\to \real^{4n-2-a}_{a+1}=\real^2_1\times \real^a\times(\real^2_1)^a\times\real^{4b}$ be the map
\begin{equation*}
f_h(t,x) = \left(\sqrt{{1}/{c}+\al(t)}\,h(\te_h(t))\,, \ti\eta(t)\,, \ti\eta\star h(t,x)\,, \bar \eta\star \vp(t,x) \right),
\end{equation*}
where  
$\te_h:I\to \real$ is the function defined by
\begin{equation*}
\te_h(t) = \int_{t_0}^t \sqrt{\frac{1}{\frac{1}{c}+\al(\tau)}\left[\rho(\tau)^2-\ep(\tau)^2 + \frac{\al'(\tau)^2}{4\left(\frac{1}{c}+\al(\tau)\right)}\right]}d\tau.
\end{equation*}
If $b=0$, we define $f_h(t,x)$ simply by omitting $\ep(t)$, $\al(t)$ and $\bar\eta\star\vp(t,x)$ in the expressions of $\te_h(t)$ and $f_h(t,x)$ above. By (\ref{cond-ep}), we have  $\rho(t)^2-\ep(t)^2>0$. Thus, in both cases $b=0$ and $b>0$, we have that $\te_h$ is well defined, smooth and increasing.

It is easy to see that $\lan f_h(t,x),f_h(t,x)\ran = -\big({1}/{c}+\al(t) \big) + |\ti\eta(t)|^2 -|\ti\eta(t)|^2 + \al(t)=  - {1}/{c},$
hence the image of $f_h$ is contained in $\hy^{4n-3-a}_a(-c)$. By using (\ref{etah}), 
\begin{eqnarray*}
f_h^*(\lan\,,\ran) &=&\left[-\frac{1}{4}\left(\frac{\al'(t)^2}{\frac{1}{c}+\al(t)}\right) + \left(\frac{1}{c}+\al(t)\right)\te'_h(t)^2 + |\ti\eta'(t)|^2\right] dt^2 \\&& + (\ti\eta\star h)^*(\lan\,,\ran) + (\bar\eta\star \vp)^*(\lan\,,\ran)\\
&=& \rho(t)^2 dt^2 + \eta_1(t)^2 dx_1^2 + \ldots + \eta_{a+b}(t)^2 dx_{a+b}^2.
\end{eqnarray*}
This implies that $f_h:M^n\to \hy^{4n-3-a}_a(-c)$ is an isometric immersion. 

Now, choose the step functions $S_1,S_2$ sufficiently large so that (\ref{cond-ep}) is satisfied and, moreover, for all $r=a+1,\ldots,a+b$, it holds 
\begin{eqnarray}\label{gransphere}
\eta_r(t)^2\psi_j(t)^2<\frac{1}{8bc}S_j(t)^2,
\end{eqnarray} 
for all $t\in I$ and $j=1,2$. By (\ref{def-a}) and (\ref{gransphere}), we obtain  $0\leq \al(t) < \frac{1}{4c}$, for all $t\in I$. 

Let $f_s:I\times \real^{n-1}\to \real_a^{4n-2-2a}=\real^2\times (\real^2_1)^a\times \real^{4b}$ be the map defined by 
\begin{equation*}
f_{s}(t,x) = \left(\sqrt{{1}/{c}-\be(t)}\,g(\te_s(t))\,, \, \ti\eta\star h\,(t,x)\,,\,\bar\eta\star \vp\,(t,x)\right),
\end{equation*}
where $\be(t)=\al(t)-|\ti\eta(t)|^2$, with $t\in I$, $g(u)=(\cos(u),\sin(u))$, with $u\in \real$, and $\te_{s}:I\to\real$ is the function given by 
\begin{equation*}
\te_s(t)= \int_{t_0}^t \sqrt{\frac{1}{\frac{1}{c}-\be(\tau)}\left[\rho(\tau)^2+|\ti\eta'(\tau)|^2-\ep(\tau)^2- \frac{\be'(\tau)^2}{4\big(\frac{1}{c}-\be(\tau)\big)}\right]}\,d\tau.
\end{equation*}
If $b=0$, we define $f_s(t,x)$  by omitting $\ep(t)$, $\al(t)$ and $\bar\eta\star\vp(t,x)$ in the expressions of $\te_s(t)$ and $f_s(t,x)$ above.
\begin{claim}\label{welldef} We can choose steps functions $S_1,S_2$, sufficiently large so that the function $\te_s$ is well defined and smooth.
\end{claim}
In fact, first we assume $b=0$. In this case, by definition, we have  $\be(t)=-|\ti\eta(t)|^2$. Hence, $\frac{1}{c}-\be(t)>|\ti\eta(t)|^2$. Moreover, $\be'(t)^2 = 4\lan \ti\eta'(t),\ti\eta(t)\ran^2\leq 4|\ti\eta'(t)|^2|\ti\eta(t)|^2$. Thus,
\begin{equation*}
\frac{\be'(t)^2}{4\big(\frac{1}{c}-\be(t)\big)}\le |\ti\eta'(t)|^2.
\end{equation*}
By (\ref{cond-ep}), we have $\rho(t)^2-\ep(t)^2>0$. Thus we conclude that $\te_s$ is well defined and smooth. Now, assume $b>0$. By (\ref{gransphere}), it holds $\frac{1}{c}-\be(t)\ge \frac{1}{c}-\al(t)>0$. Using Lemma \ref{lemmaderiv}, we have
\begin{equation}\label{deriv-a}
\al'(t) = \sum_{r=a+1}^{a+b}\left[\frac{\left(\eta_r(t)^2\psi_1(t)^2\right)'}{S_1(t)^2}+\frac{\left(\eta_r(t)^2\psi_2(t)^2\right)'}{S_2(t)^2}\right].
\end{equation}
Since $\be'(t)=\al'(t)-2\lan \ti\eta'(t),\ti\eta(t)\ran$, we obtain   
$\frac{\be'(t)^2}{4} \leq \de(t) + |\ti \eta'(t)|^2|\ti\eta(t)|^2,$ 
where $\de:I\to [0,\infty)$ is the continuous function given by
\begin{equation}\label{def-delta}
\de(t) = \left|\frac{\al'(t)^2}{4}- \al'(t)\lan\ti\eta'(t),\ti\eta(t)\ran \right|. 
\end{equation}
Using that $\frac{1}{c}-\be(t)=\frac{1}{c}-\al(t)+|\ti\eta(t)|^2>\frac{1}{2c}+|\ti\eta(t)|^2$, it holds that 
\begin{equation*}
\frac{\be'(t)^2}{4(\frac{1}{c}-\be(t))} \leq \frac{1}{\frac{1}{2c}+|\ti\eta(t)|^2}(\de(t) + |\ti\eta(t)|^2|\ti\eta'(t)|^2).
\end{equation*}
This implies that 
\begin{eqnarray*}
\rho(t)^2+|\ti\eta'(t)|^2-\ep(t)^2-\frac{\be'(t)^2}{4(\frac{1}{c}-\be(t))} &\geq& 
\rho(t)^2 - \ep(t)^2+|\ti\eta'(t)|^2\left(1-\frac{|\ti\eta(t)|^2}{\frac{1}{2c}+|\ti\eta(t)|^2}\right) \\&&  - \, \frac{\de(t)}{\frac{1}{2c}+|\ti\eta(t)|^2}\\
&=& \frac{1}{\frac{1}{2c}+|\ti\eta(t)|^2}(\Ga(t)-\de(t)),
\end{eqnarray*}
where $\Ga:I\to \real$ is the continuous function given by
\begin{equation*}
\begin{array}{l}
\Ga(t) = (\frac{1}{2c}+|\ti\eta(t)|^2)\left(\rho(t)^2 - \ep(t)^2+|\ti\eta'(t)|^2\left(1-\frac{|\ti\eta(t)|^2}{\frac{1}{2c}+|\ti\eta(t)|^2}\right)\right).
\end{array}
\end{equation*}
Since the step functions $S_1$ and $S_2$ satisfy (\ref{cond-ep}), we obtain $\Ga(t)>0$, for all $t\in I$. Furthermore, if $S_1(t)$ and  $S_2(t)$ become larger, then $\Ga(t)>0$ increases and $\de(t)$ is as smaller as we want.
So, we choose each step of $S_1$ and $S_2$ sufficiently large so that  
$\de(t)<\Ga(t)$, for all $t\in I$.  This implies that 
\begin{equation}\label{well-def-sphere}
\rho(t)^2+|\ti\eta'(t)|^2-\ep(t)^2-\frac{\be'(t)^2}{4(\frac{1}{c}-\be(t))}>0,
\end{equation}
for all $t\in I$. Claim \ref{welldef} is proved.
\\

It is easy to see that $\lan f_s(t,x),f_s(t,x)\ran = {1}/{c}$, hence the image of $f_s$ is contained in $\Sp^{4n-3-2a}_a(c)$. By using (\ref{etah}), 
\begin{eqnarray*}
f_s^*(\lan\,,\ran) &=& \left(\frac{\be'(t)^2}{4\left(\frac{1}{c}-\be(t)\right)}+\left(\frac{1}{c}-\be(t)\right)\te'_s(t)^2\right)dt^2 + (\ti\eta\star h)^*(\lan\,,\ran) + (\bar\eta\star \vp)^*(\lan\,,\ran)\\ 
&=& \rho(t)^2 dt^2 + \eta_1(t)^2 dx_1^2 + \ldots + \eta_{a+b}(t)^2 dx_{a+b}^2.
\end{eqnarray*}
This implies that $f_s:M^n \to \Sp_a^{4n-3-2a}(c)$ is an isometric immersion. Item \ref{euc-4n-2l-3} is proved.

\begin{remark}\label{non-inj} {\it  The immersions $f$, $f_h$ and $f_s$ are not injective, if $b>0$. In fact, we take $t=\mathbf{t}_{2k}$, for some integer $k$. Let  $x^1=(x^1_1,\ldots,x^1_{a+b})$ and $x^2 =(x^2_1,\ldots,x^2_{a+b})$ be vectors satisfying the following. 
\begin{enumerate}[(i)]
\item $(x^1_1,\ldots,x^1_a) = (x^2_1,\ldots,x^2_a)$;
\item $S_1(t)x^1_r=S_1(t)x^2_r + 2\pi l_{r}$, for some integer $l_r$, with $r=a+1,\ldots,a+b$ and $l_r\neq 0$ for some $r$.
\end{enumerate}
Notice that $\ti\eta\star h\,(t,x^1) = \ti\eta\star h\,(t,x^2)$, since $(x^1_1,\ldots,x^1_a) = (x^2_1,\ldots,x^2_a)$.
Since $\psi_2(t)=\psi_2(\mathbf{t}_{2k})=0$ and $(\cos(S_1(t)x^1_r),\sin(S_1(t)x^1_r)) = (\cos(S_1(t)x^2_r),\sin(S_1(t)x^2_r))$, we obtain
\begin{equation*}
\psi_j(t)(\cos(S_j(t)x^1_r),\sin(S_j(t)x^1_r)) = \psi_j(t)(\cos(S_j(t)x^2_r),\sin(S_j(t)x^2_r)).
\end{equation*}
This implies that $\bar\eta\star\vp(t,x^1)=\bar\eta\star\vp(t,x^2)$. Since the first coordinates of $f, f_h$ and $f_s$ depend only on $t$, it follows  that $f(t,x^1)=f(t,x^2)$, $f_h(t,x^1)=f_h(t,x^2)$, and $f_s(t,x^1)=f_s(t,x^2)$. Thus, the immersions  $f, f_h$ and $f_s$ are not injective.}
\end{remark}

Now we will prove Item \ref{imb-8n-6a-7}.  We will continue to assume the notations being as given in the proof of Item \ref{euc-4n-2l-3}. Let  $T_1:\real\to (0,\frac{\pi}{2})$ and $T_2:\real\to \real$ be the smooth functions 
\begin{equation} \label{T1T2}
T_1(u) = \frac{\pi}{4}\left(1+\tanh(u)\right) \ \ \mbox{ and } \ \ 
T_2(u) = \int_0^u \sqrt{1-T'_1(\tau)^2}d\tau.
\end{equation}
Note that $T_2$ is smooth since $T_1$ is analytic and $T_1'(u)=\frac{\pi}{4}\sech^2(u) \leq \frac{\pi}{4}<1$.

Consider the map 
$\hat\vp=(\vp_{11},\vp_{21},\vp_{12},\vp_{22}):I\times \real\to \real^8$, where each map $\vp_{ji}:I\times \real\to \real^2$, with $i,j=1,2$, is defined by
\begin{equation}\label{def-vp-jk}
\vp_{ji}(t,u) = \frac{\psi_j(t)}{S_j(t)}\Big(\cos\left(T_i(S_j(t)u)\right),\,\sin\left(T_i(S_j(t)u)\right)\Big).
\end{equation} 
Consider the map
\begin{equation}\label{starfunction2}
(\bar\eta\star\hat\vp)(t,x_{1},\ldots,x_{a+b}) = \big(\eta_{a+1}(t)\hat\vp(t,x_{a+1}), \ldots, \eta_{a+b}(t)\hat\vp(t,x_{a+b})\big)\in\real^{8b},
\end{equation}
with $t\in I$ and $x\in \real^{n-1}=\real^{a+b}$. 
Since $T'_1(t)^2+T_2'(t)^2= \psi_1(t)^2+\psi_2(t)^2= 1$, by using Lemma \ref{lemmaderiv}, it follows similarly as  in (\ref{etah}) that the pull-back symmetric tensor by the map $\bar\eta\star\hat\vp:I\times \real^{n-1}\to \real^{8b}$ satisfies 
\begin{equation}\label{hatvp}
(\bar\eta\star\hat\vp)^*(\lan\,,\ran) = 2\ep(t)^2\, dt^2 + \eta_{a+1}(t)^2dx_{a+1}^2+\ldots+\eta_{a+b}(t)^2dx_{a+b}^2,
\end{equation}
where $\ep:I\to [0,\infty)$ is the smooth function defined as in (\ref{ep_r}). We choose the step functions $S_1$ and $S_2$ so that  (\ref{cond-ep}) is satisfied. This implies that $\rho(t)^2-2\ep(t)^2>0$.

Let $\hat f:I\times \real^{n-1}\to \real^{8n-7-6a}_a=\real\times (\real^2_1)^a\times \real^{8b}$ be the map
\begin{equation*}
\hat f(t,x) = \left(\int_{t_0}^t\sqrt{\rho(\tau)^2+|\ti\eta'(\tau)|^2-2\ep(\tau)^2}\,d\tau\,,\, \ti\eta\star h\,(t,x)\,,\, \bar \eta\star \hat\vp\,(t,x)\right),
\end{equation*}
where $\ti\eta\star h:I\times \real^{n-1}\to (\real^2_1)^a$ is the map defined as in (\ref{starfunction}). If $b=0$, we define $\hat f(t,x)$ by simply omitting $\ep(t)$ and $\bar\eta\star\hat\vp(t,x)$ in the definition of $\hat f(t,x)$ above. By using (\ref{etah}) and (\ref{hatvp}), it is easy to conclude that  $\hat f:M^n\to \real^{8n-7-6a}_a$ is an isometric immersion. 
\begin{claim}\label{hat f inj} The immersion $\hat f$ is  injective.
\end{claim}
In fact, assume that $\hat f(t^1,x^1)=\hat f(t^2,x^2)$, for some $t^1,t^2\in I$ and $x^1,x^2\in \real^{n-1}$. We write $x^j=(x^j_1,\ldots,x^j_{a+b})$, with $j=1,2$. Using that the function 
\begin{equation*}
s(t)=\int_{t_0}^t\sqrt{\rho(\tau)^2+|\ti\eta'(\tau)|^2-2\ep(\tau)^2}\,d\tau, \ t\in I,
\end{equation*}
is increasing, we obtain  $t^1=t^2$. Since $\psi_1^2+\psi_2^2=1$, we can assume, without loss of generality, that $\psi_1(t^1)\neq 0$. Using that $\eta_i(t)>0$, for all $i=1,\ldots,a+b$ and $\hat f(t^1,x^1)=\hat f(t^1,x^2)$, we have  $h(x^1_k)=h(x^2_k)$ and $\vp_{11}(t^1,x^1_r)=\vp_{11}(t^1,x^2_r)$, for all $k=1,\ldots,a$ and $r=a+1,\ldots,a+b$. These imply that 
\begin{equation*}
\begin{array}{l}
\sinh(x^1_k)=\sinh(x^2_k)  \ \mbox{ and } \ \sin(T_1(S_1(t^1)x^1_r))=\sin(T_1(S_1(t^1)x^2_r)),
\end{array}
\end{equation*}
for all $k=1,\ldots,a$ and $r=a+1,\ldots,a+b$.
Since $S_1(t^1)>0$ and the functions $\sinh(u)$ and $\sin(T_1(u))$, with $u\in\real$, are injective, we obtain that $x^1=x^2$. Claim \ref{hat f inj} is proved.

\begin{claim}\label{hat f emb} $\hat f:M^n\to \real^{8n-7-6a}_a$ is an isometric embedding.
\end{claim}
We just need to show that the inverse map $\hat f^{-1}:\hat f(I\times \real^{n-1})\to I\times \real^{n-1}$ is continuous. In fact, let $y_m=\hat f(t_m,x^m_1,\ldots,x^m_{n-1})$ be a sequence that converges to a point $y_\infty=\hat f(t_\infty,x^\infty_1,\ldots,x^\infty_{n-1})$. Since the function $s(t)$ is the first coordinate of $\hat f(t,x)$, we obtain $\lim s(t_m)=s(t_\infty)$. This implies that $\lim t_m=t_\infty$, since $s:I\to \real$ is a diffeomorphism of $I$ onto its image $s(I)$. Since the coordinates of the map $\eta(t)=(\eta_1(t),\ldots,\eta_{n-1}(t))$ are positive and smooth, we obtain
\begin{enumerate}[(a)]
\item\label{it-c} $\lim h(x^m_k)=h(x^\infty_k)$
\item\label{it-b} $\lim \vp_{ji}(t_m,x^m_r)=\vp_{ji}(t_\infty,x^\infty_r)$,
\end{enumerate}
for all $i,j=1,2$, \, $k=1,\ldots,a$\, and \, $r=a+1,\ldots,n-1$.
It follows from \ref{it-c} that  $\lim x^m_k=x^\infty_k$, for all $k=1,\ldots,a$, since $h(u)=(\cosh(u),\sinh(u))$ and $\sinh(u)$, with $u\in \real$, is a diffeomorphism. Now, using that $\psi_1(t_\infty)^2+\psi_2(t_\infty)^2=1$, we can assume that $\psi_1(t_\infty)\neq 0$. Since $\psi_1(t_\infty)>0$, we obtain that $S_1$ is a positive constant function in a neighborhood of $t_\infty$. This implies that $S_1(t_m)=S_1(t_\infty)>0$, for  sufficiently large $m$. Since $\lim \psi_1(t_m)=\psi_1(t_\infty)>0$, we obtain from \ref{it-b} and (\ref{def-vp-jk}) that 
\begin{eqnarray*}
\lim \cos(T_1(S_1(t_\infty)x^m_r))&=& \lim \frac{S_1(t_m)}{\psi_1(t_m)}P (\vp_{11}(t_m,x^m_r)) = \frac{S_1(t_\infty)}{\psi_1(t_\infty)}P(\vp_{11}(t_\infty,x^\infty_r)) \\&=& \cos(T_1(S_1(t_\infty)x^\infty_r)),
\end{eqnarray*}
for all $r=a+1,\ldots,n-1$, where $P:\real^2\to \real$ is the projection $P(u,v)=u$. Again using that $S_1(t_\infty)> 0$ and since  $\cos(T_1(u))$ is a diffeomorphism of $\real$ onto $(0,1)$, it follows that $\lim x^m_r=x^\infty_r$, for all $r=a+1,\ldots,n-1$. We conclude that $\hat f^{-1}$ is continuous. Claim \ref{hat f emb} is proved.  
\\

Let $\hat f_h: I\times \real^{n-1}\to \real^{8n-6-5a}_{a+1}=\real^2_1\times \real^a\times (\real^2_1)^a\times \real^{8b}$ be the map
\begin{equation}\label{def hat f h}
\hat f_h(t,x) = \left(\sqrt{\frac{1}{c}+2\al(t)}\,h(\hat\te_h(t))\,,\, \ti\eta(t)  \,,\,\, \ti\eta\star h\,(t,x)\,,\, \bar \eta\star \hat\vp\,(t,x)\right),
\end{equation}
where $\al:I\to [0,\infty)$ is as defined in (\ref{def-a}) and $\hat\te_h:I\to \real$ is the function given by
\begin{equation}
\hat\te_h(t)=\int_{t_0}^t \sqrt{\frac{1}{\frac{1}{c}+2\al(\tau)}\left(\rho(\tau)^2-2\ep(\tau)^2+\frac{\al'(\tau)^2}{\frac{1}{c}+2\al(\tau)}\right)}\,d\tau.
\end{equation}
If $b=0$, we define $\hat f_h(t,x)$ simply by omitting $\al(t)$, $\ep(t)$ and $\bar\eta\star\vp(t,x)$ in the expressions of $\hat\te_h(t)$ and $\hat f_h(t,x)$ above. 

By (\ref{cond-ep}), we have $\rho(t)^2-2\ep(t)^2>0$. This implies that $\hat\te_h:I\to \real$ is well defined, smooth and increasing. 

Note that $\lan \hat f_h(t,x)\hat f_h(t,x)\ran = - \left(\frac{1}{c}+2\al(t)\right) + |\ti\eta(t)|^2 - |\ti\eta(t)|^2 + 2\al(t) = -\frac{1}{c}$. Thus the image of $\hat f_h$ is contained in $\hy^{8n-7-5a}_a(-c)$. By a standard computation, 
\begin{eqnarray*}
(\hat f_h)^*(\lan\,,\ran) &=& \left(-\frac{\al'(t)^2}{\frac{1}{c}+2\al(t)} + \left(\frac{1}{c}+2\al(t)\right)\hat\te_h'(t)^2+|\ti\eta'(t)|^2\right) dt^2 \\&& \,+\, (\ti\eta\star h)^*(\lan\,,\ran) + (\bar\eta\star \hat\vp)^*(\lan\,,\ran)\\&=& \rho(t)^2 dt^2 + \eta_1(t)^2 dx_1^2+\ldots+\eta_{a+b}(t)^2dx_{a+b}^2.
\end{eqnarray*}
Thus $\hat f_h:M^n\to \hy^{8n-7-5a}_a(-c)$ is an isometric immersion. 

\begin{claim}\label{hat f h inj} The immersion $\hat f_h$ is injective.
\end{claim}
In fact, assume that $\hat f_h(t^1,x^1)=\hat f_h(t^2,x^2)$. Using (\ref{def hat f h}), we have 
\begin{equation*}
\sqrt{\frac{1}{c}+2\al(t^1)}\,h(\hat\te_h(t^1)) = \sqrt{\frac{1}{c}+2\al(t^2)}\,h(\hat\te_h(t^2)).
\end{equation*}
Since $\lan h(u), h(u)\ran=-1$, for all $u\in \real$, we obtain $\frac{1}{c}+2\al(t^1)=\frac{1}{c}+2\al(t^2)$, hence $h(\hat\te_h(t^1))=h(\hat\te_h(t^2))$. This implies that $t^1=t^2$, since the function $\sinh(\hat\te_h(t))$, with $t\in I$, is increasing.  The argument to show that $x^1=x^2$ is similar the one as given in Claim \ref{hat f inj}. Claim \ref{hat f h inj} is proved.

\begin{claim}\label{hat f h emb} $\hat f_h:M^n\to \hy^{8n-7-5a}_a(-c)$ is an isometric embedding.
\end{claim}
In fact,  we just need to prove that the inverse map $(\hat f_h)^{-1}: \hat f_h(I\times \real^{n-1})\to I\times \real^{n-1}$ is continuous. Let $y_m=\hat f_h(t_m,x^m)$ be a sequence that converges to a point $y_\infty=\hat f_h(t_\infty,x^\infty)$. Using (\ref{def hat f h}), we obtain 
\begin{equation*}
\lim \sqrt{\frac{1}{c}+2\al(t_m)}\, h(\hat\te_h(t_m))=\sqrt{\frac{1}{c}+2\al(t_\infty)}\, h(\hat\te_h(t_\infty)).
\end{equation*}
Using $\lan h(u),h(u)\ran=-1$, we obtain $\lim (\frac{1}{c}+2\al(t_m))=\frac{1}{c}+2\al(t_\infty)>0$, hence $\lim h(\hat\te_h(t_m))=h(\hat\te_h(t_\infty))$. This implies that $\lim \sinh(T_1(\hat\te_h(t_m)))=\sinh(T_1(\hat\te_h(t_\infty)))$. Using that $\sinh(T_1(\hat\te_h(t)))$ is a diffeomorphism of $I$ onto its image, it follows that $\lim t_m=t_\infty$. The argument to show that $\lim x^m=x^\infty$ is also similar to Claim \ref{hat f emb}. Thus Claim \ref{hat f h emb} is proved. 
\\

Now let $\hat f_s:I\times \real^{n-1}\to \real^{8n-4-6a}_a=\real^4\times (\real^2_1)^a\times \real^{8b}$ be the map
\begin{equation}\label{def-hat f s}
\hat f_s(t,x) = \left(\sqrt{\frac{1}{2c}+\frac{|\ti\eta(t)|^2}{2}-\al(t)}\, \m C(\hat\te_s(t))\,,\,\,\ti\eta\star h\,(t,x)\,,\,\bar\eta\star \hat\vp\,(t,x)\right),
\end{equation}
where $\m C(u)=(\cos(T_1(u)),\sin(T_1(u)), \cos(T_2(u)),\sin(T_2(u)))$, with $u\in \real$. Further, the function $\al:I\to   [0,\infty)$ is given as in (\ref{def-a}),
and $\hat\te_s:I\to \real$ is defined by
\begin{equation*}
\hat\te_s(t)=\int_{t_0}^t \sqrt{\frac{G(\tau)}{\frac{1}{2c}+\frac{|\ti\eta(\tau)|^2}{2} - \al(\tau)} }\,d\tau,
\end{equation*}
where $G:I\to \real$ is the function $$G(t)=\rho(t)^2 +|\ti\eta'(t)|^2 - 2\ep(t)^2 - 2\frac{\left(\left(-\al(t)+\frac{1}{2}|\ti\eta(t)|^2\right)'\right)^2}{4(\frac{1}{2c}+\frac{|\ti\eta(t)|^2}{2}-\al(t))}.$$
If $b=0$, we define $\hat f_s$ by simply omitting $\al(t)$ and $\bar\eta\star\hat\vp(t,x)$ in the definitions of $\hat\te_s(t)$ and $\hat f_s(t,x)$ above.

We claim that we can choose the step functions $S_1$ and $S_2$ sufficiently large so that $\hat\te_s$ is well defined and smooth. By (\ref{gransphere}), we already have $\frac{1}{2c}-\al(t)>0$. Furthermore, by a simple computation, 
\begin{eqnarray*}
G(t) &=& \rho(t)^2 + |\ti\eta'(t)|^2 - \frac{\lan\ti\eta'(t),\ti\eta(t)\ran^2}{\frac{1}{c}+|\ti\eta(t)|^2-2\al(t)} - \De(t)\\
&\ge & \rho(t)^2 + |\ti\eta'(t)|^2 \left[1-\frac{|\ti\eta(t)|^2}{\frac{1}{c}+|\ti\eta(t)|^2-2\al(t)}\right] - \De(t),\end{eqnarray*}
where $\De(t)=2\ep(t)^2 + \frac{\al'(t)^2 - 2\al'(t)\lan \ti\eta'(t),\ti\eta(t)\ran}{\frac{1}{c}+|\ti\eta(t)|^2-2\al(t)}$.
Note that we can take $\ep(t)$, $\al(t)$ and $\al'(t)$ as smaller as we want if $S_1(t)$ and $S_2(t)$ become larger. Thus, we can choose the step functions $S_1$ and $S_2$ sufficiently large so that $\De(t)<\rho(t)^2$. This implies that $\hat \te_s (t)$ is well defined, smooth and increasing.

Note that $\lan \hat f_s(t,x),\hat f_s(t,x)\ran = \frac{1}{c}$ since $|\m C(u)|^2=2$,  $|\bar\eta\star \hat\vp\,(t,x)|^2=2\al(t)$ and  $\lan \ti\eta\star h\,(t,x)\,,\,\ti\eta\star h\,(t,x) \ran=-|\ti\eta(t)|^2$. Thus the image $\hat f_s(I\times \real^{n-1})\subset \Sp^{8n-5-6a}_a(c)$. 

By a direct computation, we show that
\begin{eqnarray*}
\hat f_s^*(\lan\,,\ran) &=& \left(\frac{\left(\left(-\al(t)+\frac{1}{2}|\ti\eta(t)|^2\right)'\right)^2}{\frac{1}{c}+|\ti\eta(t)|^2-2\al(t)} + \left(\frac{1}{2c}+\frac{|\ti\eta(t)|^2}{2}-\al(t)\right)\hat\te_s'(t)^2\right)dt^2\\&& + (\ti\eta\star h)^*(\lan\,,\ran) + (\bar\eta\star \hat\vp)^*(\lan\,,\ran)\\&=& \rho(t)^2 dt^2 + \eta_1(t)^2 dx_1^2 +\ldots+\eta_{a+b}(t)^2 dx_{a+b}^2.
\end{eqnarray*}
This implies that $\hat f_s:M^n\to \Sp^{8n-5-6a}_a(c)$ is an isometric immersion.
\begin{claim}\label{hat f s inj} The immersion $\hat f_s$ is injective.
\end{claim}
In fact, assume that $\hat f_s(t^1,x^1)=\hat f_s(t^2,x^2)$, for some $t^1,t^2\in I$ and $x^1,x^2\in \real^{n-1}$. Using (\ref{def-hat f s}), $$\sqrt{\frac{1}{c}+|\ti\eta(t^1)|^2-2\al(t^1)}\, \m C(\hat\te_s(t^1))=\sqrt{\frac{1}{c}+|\ti\eta(t^2)|^2-2\al(t^2)}\, \m C(\hat\te_s(t^2)).$$
Since $|\m C(\hat\te_s(t))|^2=2$, for all $t\in I$, we obtain  $|\ti\eta(t^1)|^2-2\al(t^1)=|\ti\eta(t^2)|^2-2\al(t^2)$, hence $\m C(\hat\te_s(t^1))=\m C(\hat\te_s(t^2))$. This implies that $t^1=t^2$, since  $\sin(T_1(\hat\te_s(t)))$, with $t\in I$, is increasing. The argument to show that $x^1=x^2$ is similar to that one given in Claim \ref{hat f inj}.
\\

\begin{claim}\label{hat f s emb} $\hat f_s:M^n\to \Sp^{8n-5-6a}_a(c)$ is an isometric embedding.
\end{claim}
In fact, we just need to prove that the inverse map $(\hat f_s)^{-1}: \hat f_s(I\times \real^{n-1})\to I\times \real^{n-1}$ is continuous. Let $y_m=\hat f_s(t_m,x^m)$ be a sequence that converges to a point $y_\infty=f(t_\infty,x^\infty)$. Using (\ref{def-hat f s}),  $$\lim \sqrt{\frac{1}{c}+|\ti\eta(t_m)|^2-2\al(t_m)}\, \m C(\hat\te_s(t_m))=\sqrt{\frac{1}{c}+|\ti\eta(t_\infty)|^2-2\al(t_\infty)}\, \m C(\hat\te_s(t_\infty)).$$
Since $|\m C(\hat\te_s(t))|^2=2$, for all $t\in I$, we have $\lim (|\ti\eta(t_m)|^2- 2\al(t_m))=|\ti\eta(t_\infty)|^2- 2\al(t_\infty)>0$, hence  $\lim \m C(\hat\te_s(t_m))=\m C(\hat\te_s(t_\infty))$. This implies that $\lim \sin(T_1(\hat\te_s(t_m)))=\sin(T_1(\hat\te_s(t_\infty)))$. Using that $\sin(T_1(\hat\te_s(t)))$ is a diffeomorphism of $I$ onto its image, it follows that $\lim t_m=t_\infty$. The argument to show that $\lim x^m=x^\infty$ is similar to Claim \ref{hat f emb}. Thus, Claim \ref{hat f s emb} is proved. Theorem \ref{mth} is proved.

\begin{small}
\vspace{1cm}
\begin{tabular}{l}
Heudson Mirandola\\
Universidade Federal do Rio de Janeiro\\
Instituto de Matem\'{a}tica\\
21945-970 Rio de Janeiro-RJ\\ Brazil\\
\verb+mirandola@im.ufrj.br+\\
\end{tabular}\
\vspace{0.3cm}
\begin{tabular}{l}
Feliciano Vit\'orio\\
Universidade Federal de Alagoas\\
Instituto de Matem\'{a}tica\\
57072-900 Macei\'o-AL\\ Brazil\\
\verb+feliciano@pos.mat.ufal.br+\\
\end{tabular}\
\end{small}

\end{document}